\documentclass{amsart}

\usepackage{enumerate}
\usepackage{amsmath}%
\usepackage{amsfonts}%
\usepackage{amssymb}%
\usepackage{graphicx}
\usepackage{mathrsfs}
\usepackage{hyperref}
\usepackage{fullpage}
\newtheorem{theorem}{Theorem}
\theoremstyle{plain}

\newtheorem{claim}[theorem]{Claim}

\newtheorem{conjecture}[theorem]{Conjecture}

\newtheorem{fact}[theorem]{Fact}
\newtheorem{lemma}[theorem]{Lemma}

\newtheorem{problem}[theorem]{Problem}

\numberwithin{equation}{section}
\numberwithin{theorem}{section}
\numberwithin{case}{section}

\numberwithin{subcase}{case}

\def\cP{\mathcal{P}}

\def \a{\alpha}
\def \e{\epsilon}
\def \r{\gamma}

\def \bfi{\mathbf{i}}
\def \bfu{\mathbf{u}}
\def \bfv{\mathbf{v}}

\begin{document}
\title{On Perfect Matchings and tilings in uniform Hypergraphs}
\author{Jie Han}
\address{Instituto de Matem\'{a}tica e Estat\'{\i}stica, Universidade de S\~{a}o Paulo, Rua do Mat\~{a}o 1010, 05508-090, S\~{a}o Paulo, Brazil}
\email[Jie Han]{jhan@ime.usp.br}
\thanks{The author was supported by FAPESP (2014/18641-5, 2013/03447-6).}
\date{\today}

\begin{abstract}
In this paper we study some variants of Dirac-type problems in hypergraphs.
First, we show that for $k\ge 3$, if $H$ is a $k$-graph on $n\in k\mathbb N$ vertices with independence number at most $n/p$ and minimum codegree at least $(1/p+o(1))n$, where $p$ is the smallest prime factor of $k$, then $H$ contains a perfect matching.
Second, we show that if $H$ is a $3$-graph on $n\in 3\mathbb N$ vertices which does not contain any induced copy of $K_4^-$ (the unique $3$-graph with $4$ vertices and $3$ edges) and has minimum codegree at least $(1/3+o(1)))n$, then $H$ contains a perfect matching.
Moreover, if we allow the matching to miss at most $3$ vertices, then the minimum degree condition can be reduced to $(1/6+o(1)))n$.
Third, we show that if $H$ is a $3$-graph on $n\in 4\mathbb N$ vertices which does not contain any induced copy of $K_4^-$ and has minimum codegree at least $(1/8+o(1)))n$, then $H$ contains a perfect $Y$-tiling, where $Y$ represents the unique $3$-graph with $4$ vertices and $2$ edges.
We also provide the examples showing that our minimum codegree conditions are asymptotically best possible.
Our main tool for finding the perfect matching is a characterization theorem in~\cite{Han14_poly} that characterizes the $k$-graphs with minimum codegree at least $n/k$ which contain a perfect matching.
\end{abstract}

\maketitle

\section{Introduction}
Matchings are fundamental objects in graphs and hypergraphs and the existence of matchings or perfect matchings has been an active area of research for many years. 
Recently, the problem of finding perfect matchings and tilings in uniform hypergraphs has received much attention from various researchers.
One main motivation to study the hypergraph matching theory is its steep contrast to the graph case: the decision problem for perfect matchings in hypergraphs is NP-complete (Karp~\cite{Karp}), whilst the corresponding problem for ordinary graphs can be solved in polynomial time (Edmonds's algorithm~\cite{Edmonds}).
Since the general decision problem is intractable provided that $\text{P}\neq \text{NP}$, it is natural to seek sufficient conditions for the existence of perfect matchings.
One well-studied class of such conditions are minimum degree conditions.
To discuss the results let us mention the following standard definitions.

Given $k\ge 2$, a $k$-uniform hypergraph (in short, \emph{$k$-graph}) $H=(V(H), E(H))$ consists of a vertex set $V(H)$ and an edge set $E(H)\subseteq \binom{V(H)}{k}$, where every edge is a $k$-element subset of $V(H)$. 
A \emph{subgraph} $H'$ of $H$ is a $k$-graph such that $V(H')\subseteq V(H)$ and $E(H')\subseteq E(H)$.
(Here we prefer the term `subgraph' instead of `subhypergraph' or `sub-$k$-graph'.)
An \emph{induced subgraph} $H'$ of $H$ is a $k$-graph such that $V(H')\subseteq V(H)$ and $E(H')=\{e\in E(H): e\subseteq V(H')\}$.
We write $H[A]$ for the induced subgraph of $H$ on the vertex set $A$.
A \emph{matching} in $H$ is a collection of vertex-disjoint edges of $H$. A \emph{perfect matching} in $H$ is a matching that covers $V(H)$. Clearly a perfect matching in $H$ exists only if $k$ divides $|V(H)|$. 
Given a $k$-graph $H$ with a set $S$ of $d$ vertices (where $1 \le d \le k-1$) we define $N_{H} (S)$ to be the collection of $(k-d)$-sets $T$ such that $S\cup T\in E(H)$, and let $\deg_H(S):=|N_H(S)|$ (the subscript $H$ is omitted if it is clear from the context). The \emph{minimum $d$-degree $\delta _{d} (H)$} of $H$ is the minimum of $\deg_{H} (S)$ over all $d$-vertex sets $S$ in $H$.  
We refer to $\delta _{k-1} (H)$ as the \emph{minimum codegree} of $H$.
The \emph{independence number of $H$}, denoted by $\alpha(H)$, is the maximum size of a set in $H$ that spans no edge.
Throughout this paper, $o(1)$ represents a function that tends to $0$ as $n$ goes to infinity.

In general, Dirac-type problems ask for \emph{the minimum $d$-degree threshold}, namely, the smallest possible minimum $d$-degree condition that ensures certain (spanning) subgraphs in a $k$-graph (such a threshold always exists because subgraph containment properties are monotone).
For example, given $d, k\in \mathbb N$ such that $d\le k-1$, one main question is to determine the minimum $d$-degree threshold that ensures a perfect matching in a $k$-graph.
For the case $d=k-1$ and sufficiently large $n\in k\mathbb{N}$, R\"odl, Ruci\'nski and Szemer\'edi~\cite{RRS09} showed that the minimum codegree threshold is $n/2-k+C$, where $C\in\{3/2, 2, 5/2, 3\}$ depends on the values of $n$ and $k$.
For general $d$, it is conjectured in~\cite{HPS, KuOs-survey} that the threshold is $\max \{ 1/2+o(1), 1-\left (1-{1}/{k}\right)^{k-d}+o(1) \}\binom{n-d}{k-d}$.
This conjecture has been confirmed for all $0.42 k\le d\le k-1$ and for $1\le k-d\le 4$ and for $(k, d)\in \{(12,5), (17,7)\}$ (see \cite{AFHRRS, CzKa, HPS, JieNote, Khan1, Khan2, KO06mat, KOT, MaRu, Pik, RRS06mat, RRS09, TrZh12, TrZh13, TrZh15}). 
Moreover, the \emph{exact} thresholds for large $n$ are obtained for these cases except $(k, d)\in \{(5,1), (6,2)\}$ in \cite{CzKa, JieNote, Khan1, Khan2, KOT, TrZh12, TrZh13, TrZh15}. 
For more results we refer the reader to the recent surveys of R\"odl and Ruci\'nski~\cite{RR}, and Zhao~\cite{zsurvey}.

The main focus of this paper is to study two variants of Dirac-type problems.
The first one is to investigate the perfect matching problems on $k$-graphs which have small independence number.
Our motivation is the recent work of Balogh--Molla--Sharifzadeh~\cite{BMS} and Balogh--McDowell--Molla--Mycroft~\cite{BMMM}, who studied triangle-factors (a spanning subgraph consisting of vertex-disjoint copies of triangles) under the minimum degree condition in graphs with sublinear independence number.
More precisely, a classical result of Corr\'adi and Hajnal~\cite{CoHa} shows that any $n$-vertex graph $G$ with minimum degree $\delta(G)\ge 2n/3$ contains a triangle-factor, and the minimum degree condition is best possible.
Balogh--Molla--Sharifzadeh~\cite{BMS} showed that if in addition $\alpha(G)=o(n)$, then one can weaken the minimum degree condition to $(1/2+o(1))n$.
More recently, Balogh--McDowell--Molla--Mycroft~\cite{BMMM} showed that any $n$-vertex graph $G$ with $\delta(G)\ge (1/3+o(1))n$ and $\alpha(G)=o(n)$ contains a triangle-tiling that leaves a constant number of vertices uncovered.

The second variant is to study the $k$-graphs with forbidden induced subgraphs.
For graphs, this is a broadly studied and fruitful area.
A classical result of Sumner~\cite{Sumner} and, independently, Las Vergnas~\cite{LasVergnas} says that if a connected graph of even order contains no induced copy of $K_{1,3}$, then it has a perfect matching.
In general there is a huge literature of research on finding (spanning) subgraphs in graphs which do not contain any induced copy of a certain graph (or a graph from a certain class of graphs).
Regarding hypergraphs, Razborov~\cite{Raz10} showed that every $3$-graph on $n$ vertices that does not contain complete subgraphs on $4$-vertices, and in which no $4$ vertices span exactly one edge, must have at most $(5/9+o(1))\binom n3$ edges.
Later Pikhurko~\cite{Pik11} determined this extremal function exactly and described all extremal hypergraphs of large order.

From the extremal perspective, the idea behind these types of problems is that if we put certain structural restrictions on the graph or $k$-graph such that the extremal configuration(s) for the original problem are forbidden, then one might expect to reduce the minimum degree (or the edge density) assumption significantly.
In this paper we consider Dirac-type problems under the minimum codegree condition together with either the independence number condition or the condition that the $k$-graph contains no induced copy of a certain $k$-graph.

\subsection{$k$-graphs with small independence number}
First we are interested in $n$-vertex $k$-graphs $H$ with small independence number.
Trivially, any $k$-graph $H$ contains a matching of size at least $(n - \a(H))/k$.
So if $\alpha(H)=o(n)$, then $H$ contains a matching of size at least $(1-o(1))n/k$.
We show that by requiring a minimum codegree condition, we can guarantee a matching that leaves at most $k$ vertices uncovered.

\begin{theorem}\label{thm:NPM1}
For any $k\ge 3$ and $\r>0$, there exists $\e>0$ such that the following holds for sufficiently large integer $n$.
Suppose $H$ is an $n$-vertex $k$-graph such that $\a(H)\le \e n$ and $\delta_{k-1}(H) \ge \r n$, then $H$ contains a matching that covers all but at most $k$ vertices.
\end{theorem}

Recall that the minimum codegree threshold for perfect matchings is about $n/2$.
We next show that if the $k$-graph $H$ has small independence number then the threshold can be reduced to $(1/p+o(1))n$, where $p$ is the smallest prime factor of $k$.
In fact, we only need $\a(H) \le n/p$ instead of $\a(H) = o(n)$.

\begin{theorem}\label{thm:main1}
For any $k\ge 3$ and $\r>0$, the following holds for sufficiently large integer $n\in k\mathbb{N}$.
Let $p$ be the smallest prime factor of $k$.
Suppose $H$ is an $n$-vertex $k$-graph such that $\a(H)\le n/p$ and $\delta_{k-1}(H) \ge (1/p+\r)n$, then $H$ contains a perfect matching.
\end{theorem}

The minimum codegree condition in Theorem~\ref{thm:main1} is tight up to the error term $\r$ (see Section 4).

\subsection{$3$-graphs without induced $K_4^-$}
We say a $k$-graph $H$ is \emph{induced $F$-free}, if it does not contain any copy of $F$ as an induced subgraph.
Let $K_4^-$ be the unique $3$-graph with $4$ vertices and $3$ edges.
Here we are interested in $3$-graphs which are induced $K_4^-$-free.
We show that if the $k$-graph $H$ is induced $K_4^-$-free then the minimum codegree threshold for perfect matchings can be reduced significantly.

\begin{theorem}\label{thm:main2}
For any $\r>0$, the following holds for sufficiently large integer $n\in 3\mathbb{N}$.
Suppose $H$ is an $n$-vertex $3$-graph which is induced $K_4^-$-free and $\delta_{2}(H) \ge (1/3+\r)n$, then $H$ contains a perfect matching.
\end{theorem}

Next we show that if we only target on a matching that is almost perfect,  then the minimum codegree threshold can be further reduced.

\begin{theorem}\label{thm:NPM2}
For any $\r>0$, the following holds for sufficiently large integer $n$.
Suppose $H$ is an $n$-vertex $3$-graph such that $H$ is induced $K_4^-$-free and $\delta_{2}(H) \ge (1/6+\r)n$, then $H$ contains a matching that covers all but at most $3$ vertices.
\end{theorem}

In general, given $k$-graphs $F$ and $H$, a \emph{perfect $F$-tiling} in a $k$-graph $H$ is a spanning subgraph of $H$ that consists of vertex-disjoint copies of $F$.
Clearly a perfect $F$-tiling exists only if $|V(H)|$ is divisible by $|V(F)|$.
Let $Y$ be the unique $3$-graph with $4$ vertices and $2$ edges.
Note that the $3$-graph has been denoted as $F$, $C_4^3$ or $K_4^3-2e$ by other researchers.
Perfect $Y$-tilings in $3$-graphs were first studied by K\"uhn and Osthus \cite{KO}, who showed that minimum codegree $(1+o(1))n/4$ ensures the existence of a perfect $Y$-tiling in an $n$-vertex $3$-graph.
Subsequently Czygrinow, DeBiasio and Nagle~\cite{CDN} determined the exact minimum codegree threshold for large $n$.
Recently the author and Zhao~\cite{HZ3} and independently Czygrinow~\cite{Czy14} determined the exact vertex degree threshold for perfect $Y$-tilings for sufficiently large $n$.
In this paper we show a result on perfect $Y$-tilings in induced $K_4^-$-free $3$-graphs.

\begin{theorem}\label{thm:main3}
For any $\r>0$, the following holds for sufficiently large integer $n\in 4\mathbb{N}$.
Suppose $H$ is an $n$-vertex $3$-graph such that $H$ is induced $K_4^-$-free and $\delta_{2}(H) \ge (1/8+\r)n$, then $H$ contains a perfect $Y$-tiling.
\end{theorem}

Again the minimum codegree conditions in Theorems~\ref{thm:main2} --~\ref{thm:main3} are asymptotically sharp (see Section 4).

Our main tool for finding perfect matchings in Theorems~\ref{thm:main1} and~\ref{thm:main2} is Theorem~\ref{lem:PL} from~\cite{Han14_poly} that characterizes the $n$-vertex $k$-graphs with minimum codegree at least $n/k$ and not so large independence number which contain a perfect matching.
The result was used in~\cite{Han14_poly} to prove that the decision problem for the existence of a perfect matching in $n$-vertex $k$-graphs with minimum codegree $n/k$ is can be solved in polynomial-time, which improved a result of Keevash--Knox--Mycroft~\cite{KKM13}.
Its proof uses a variant of the absorbing method, named `lattice-based absorbing method' developed by the author, which has been applied to solve other matching and tiling problems~\cite{GaHa, GHZ, Han15_mat, Jieconf, HLTZ_K4, HT, HZZ_mat, HZZ_tiling}. 



\section{Proofs of Theorems~\ref{thm:NPM1}--\ref{thm:main2}}

\subsection{Proofs of Theorem~\ref{thm:NPM1}}

The \emph{absorbing method} initiated by R\"odl, Ruci\'nski and Szemer\'edi~\cite{RRS09} has proven to be an efficient tool for finding spanning structures in graphs and hypergraphs.
We start with a simple absorbing lemma for near perfect matchings from~\cite{RRS09}.

Let $H$ be a $k$-graph.
Given a set $S$ of $k+1$ vertices of $H$, we call an edge $e\in E(H)$ disjoint from $S$ \emph{$S$-absorbing} if there are two disjoint edges $e_1$ and $e_2$ in $E(H)$ such that $|e_1\cap S|=k - 1$, $|e_1\cap e| = 1$, $|e_2\cap S|=2$, and $|e_2\cap e|=k-2$. 
Let us explain how such absorbing works. Let $S$ be a $(k+1)$-set and $M$ be a matching, where $V(M)\cap S=\emptyset$, which contains an $S$-absorbing edge $e$. Then $M$ can ``absorb'' $S$ by replacing $e$ in $M$ by $e_1$ and $e_2$ (one vertex of $e$ becomes uncovered).
We use the absorbing lemma~\cite[Fact 2.3]{RRS09} with some quantitative changes as follows. 
We include a proof for completeness.

Throughout this paper, $x\ll y$ means that for any $y> 0$ there exists $x_0> 0$ such that for any $x< x_0$ the subsequent statement holds. Hierarchy of more constants are defined similarly.
Moreover when we include $1/n$ in the hierarchy, we implicitly mean that $n$ is an integer.

\begin{lemma}\label{lem:abs2}
Suppose $k\ge 3$ and $0<1/n \ll \beta \ll \r, 1/k$.
Let $H$ be an $n$-vertex $k$-graph with $\delta_{k-1}(H)\ge \r n$, then there exists a matching $M'$ in $H$ of size $|M'|\le \beta n$ such that for every $(k+1)$-tuple $S$ of vertices of $H$, the number of $S$-absorbing edges in $M'$ is at least $\beta^2 n$.
\end{lemma}

We need the following fact~\cite[Fact 2.2]{RRS09} in the proof of Lemma~\ref{lem:abs2}.

\begin{fact}\cite{RRS09}\label{fact0}
Suppose $0<1/n \ll \r, 1/k$.
Let $H$ be an $n$-vertex $k$-graph with $\delta_{k-1}(H)\ge \r n$, then for every set $S$ of $k+1$ vertices, $H$ contains at least $\r^3 n^k/(2k!)$ $S$-absorbing edges. 
\end{fact}

\begin{proof}[Proof of Lemma~\ref{lem:abs2}]
Select a random subset $M$ of $E(H)$, where each edge is chosen independently with probability $p=\beta n^{1-k}$.
Thus the expected size of $M$ is at most $\binom nk p \le n^k p/k! \le \beta n/6$, and the expected number of intersecting pairs of edges in $M$ is at most $n^{2k-1}p^2 = \beta^2 n$.
Hence, by Markov's inequality, $M$ contains at most $\beta n $ edges and there are at most $3\beta^2 n$ intersecting pairs of edges in $M$ with probability at least $1-1/6-1/3=1/2$.

For every $(k+1)$-tuple of vertices $S$, let $X_S$ be the number of $S$-absorbing edges in $M$.
Then by Fact~\ref{fact0}, $\mathbb{E}[X_S] \ge p \r^3 n^k/(2k!) \ge 8\beta^2 n$, because $\beta \ll \r$.
By Chernoff's bound (see~\cite[Theorem 2.1]{JLR}), we have
\[
\mathbb{P}\left( X_S\le \frac12 \mathbb{E}[X_S]\right) \le \exp \left( - \frac18 \mathbb{E}[X_S]\right)\le \exp(-\beta n^2) \le n^{-k-2},
\]
as $1/n\ll \beta$.
Thus with probability at least $1-1/n$, for every set $S$ of $k+1$ vertices, $X_S\ge 4\beta^2 n$.

Since $1/2+1/n<1$, there exists a choice of $M$ such that $|M|\le \beta n$, $M$ contains at most $3\beta^2 n$ intersecting pairs of edges, and for every set $S$ of $k+1$ vertices, $M$ contains at least $4\beta^2 n$ $S$-absorbing edges.
We obtain $M'$ by deleting one edge from each intersecting pair of edges in $M$.
Thus $M'$ is a matching in $H$, $|M'|\le \beta n$, and for every set $S$ of $k+1$ vertices, $M$ contains at least $\beta^2 n$ $S$-absorbing edges.
\end{proof}

Theorem~\ref{thm:NPM1} can be quickly derived from Lemma~\ref{lem:abs2}.

\begin{proof}[Proof of Theorem~\ref{thm:NPM1}]
Given an integer $k\ge 3$ and a real $\r>0$, suppose $0<1/n \ll \epsilon\ll \beta \ll \r, 1/k$.
Then let $H$ be an $n$-vertex $k$-graph such that $\alpha(H)\le \e n$ and $\delta_{k-1}(H)\ge \r n$.
By Lemma~\ref{lem:abs2}, $H$ contains a matching $M'$ of size at most $\beta n$ such that for every $(k+1)$-tuple $S$ of vertices of $H$, the number of $S$-absorbing edges in $M'$ is at least $\beta^2 n$.
Let $H'=H\setminus V(M')$ and note that $\alpha(H')\le \e n$.
Thus there is a matching $M''$ in $H'$ that covers all but a set $X$ of at most $\e n \le \beta^2 n$ vertices.
Since for every $(k+1)$-tuple $S$ of vertices of $H$, the number of $S$-absorbing edges in $M'$ is at least $\beta^2 n$, we can repeatedly absorb the leftover vertices (at most $\beta^2 n/k$ times, each time the number of leftover vertices is reduced by $k$) until the number of unmatched vertices is at most $k$.
Let $\tilde{M}$ denote the matching on $V(M')\cup X$ that leaves at most $k$ vertices unmatched.
Then $\tilde{M}\cup M''$ is the desired matching in $H$.
\end{proof}

\subsection{Proofs of Theorems~\ref{thm:main1} and \ref{thm:main2}}
The proofs of Theorems~\ref{thm:main1} and~\ref{thm:main2} are more involved. 
We first motivate the notion of the lattice structures needed for our proof.
For a $k$-graph $H$, we will find a partition $\cP$ of $V(H)$ and then study the so-called robust edge-lattice with respect to this partition.
In particular, every partition $\cP$ considered in this paper includes an order on its parts.

The following definitions are introduced by Keevash and Mycroft \cite{KM1}. 
Fix an integer $d>0$, let $H=(V,E)$ be a $k$-graph and let $\cP=\{V_1,\dots, V_d\}$ be a partition of $V$.
The \emph{index vector} $\mathbf{i}_{\cP}(S)\in \mathbb{Z}^d$ of a subset $S\subset V$ with respect to $\cP$ is the vector whose coordinates are the sizes of the intersections of $S$ with each part of $\cP$, i.e., $\mathbf{i}_{\cP}(S)_{V_i}=|S\cap V_i|$ for $i\in [d]$.
We call a vector $\mathbf{i}\in \mathbb{Z}^d$ an \emph{$s$-vector} if all its coordinates are nonnegative and their sum equals $s$. 
Given $\mu>0$, a $k$-vector $\mathbf{v}$ is called a $\mu$\emph{-robust edge-vector} if at least $\mu |V|^k$ edges $e\in E$ satisfy $\mathbf{i}_\cP(e)=\mathbf{v}$.
Let $I_{\cP}^{\mu}(H)$ be the set of all $\mu$-robust edge-vectors and let $L_{\cP}^{\mu}(H)$ be the lattice (additive subgroup) generated by the vectors of $I_{\cP}^{\mu}(H)$.
For $j\in [d]$, let $\mathbf{u}_j\in \mathbb{Z}^d$ be the $j$-th \emph{unit vector}, namely, $\mathbf{u}_j$ has 1 on the $j$-th coordinate and 0 on other coordinates.
A \emph{transferral} is the vector $\bfu_i - \bfu_j$ for some $i\neq j$.
A lattice $L\subseteq \mathbb{Z}^d$ is called \emph{transferral-free} if for any $i\neq j$, $\bfu_i- \bfu_j\notin L$.
We say that $L$ is \emph{full} if it is transferral-free and for every $(k-1)$-vector $\bfv$ there is some $i\in [d]$ such that $\bfv+\bfu_i\in L$.
%
Let $\mu >0$, then we say $(\cP, L_{\cP}^\mu (H))$ is \emph{full} if $|\cP|\le k$ and $L_{\cP}^\mu(H)$ is full.

Moreover, we need the following definition and result from~\cite{KKM13}.
Let $L\subseteq \mathbb{Z}^d$ be a lattice generated by a set of $k$-vectors in $\mathbb{Z}^d$ and let $L_{\max}^d$ be the lattice generated by the set of all $k$-vectors in $\mathbb{Z}^d$.
Then the \emph{coset group} of $L$ is $G=G(L)=L_{\max}^d/L$. 

\begin{lemma}\cite[Lemma 6.4]{KKM13}\label{lem:group}
If $k\ge 3$ and $L\subseteq \mathbb{Z}^d$ is a full lattice generated by a set of $k$-vectors in $\mathbb{Z}^d$, then $|G(L)|=d$.
\end{lemma}

We also use the reachability arguments introduced by Lo and Markstr\"om \cite{LM2, LM1}.
We say that two vertices $u$ and $v$ are \emph{$(\beta, i)$-reachable} in $H$ if there are at least $\beta n^{i k-1}$ $(i k-1)$-sets $S$ such that both $H[S\cup \{u\}]$ and $H[S\cup \{v\}]$ have perfect matchings. 
We say that a vertex set $U$ is \emph{$(\beta, i)$-closed in $H$} if any two vertices $u,v\in U$ are $(\beta, i)$-reachable in $H$.
For two partitions $\cP, \cP'$ of a set $V$, we say that $\cP$ \emph{refines} $\cP'$ if every vertex class of $\cP$ is a subset of some vertex class of $\cP'$. 

We combine~\cite[Lemma 6.9]{KKM13},~\cite[Lemma 2.5]{Han14_poly} and~\cite[Theorem 2.8]{Han14_poly} in the following theorem.
We say that $H$ is \emph{$\r$-extremal} if $V(H)$ contains an independent subset of size at least $(1-\r)\frac{k-1}k n$. 
A (possibly empty) matching $M$ of size at most $|\cP|-1$ is a \emph{solution} for $(\cP, L)$ (in $H$) if $\bfi_{\cP}(V\setminus V(M))\in L$; we say that $(\cP, L)$ is \emph{soluble} if it has a solution.

\begin{theorem}\cite{Han14_poly}\label{lem:PL}
Given an integer $k\ge 3$, suppose that
$1/n \ll \beta, \mu_0 \ll \r' \ll 1/k$.
Then for each $k$-graph $H$ on $n$ vertices with $\delta_{k-1}(H)\ge n/k-\r' n$, there exist $\mu_0^2 \le \mu\le \mu_0$ and partitions $\cP=\{V_1, \dots, V_d \}$ and $\cP'=\{V_1', \dots, V_{d'}'\}$ of $V(H)$ satisfying the following properties:
\begin{enumerate}[\emph{(}i\emph{)}]
\item $\cP'$ refines $\cP$ and $(\cP, L_{\cP}^{\mu}(H))$ is a full pair,
\item each partition set of $\cP$ or $\cP'$ has size at least $n/k-2\r' n$,
\item for each $D\subseteq V(H)$ such that $\bfi_{\cP}(D)\in L_{\cP}^\mu(H)$, we have $\bfi_{\cP'}(D)\in L_{\cP'}^\mu(H)$,
\item for each $i\in [d]$, $V_i'$ is $(\beta, 2^{k-1})$-closed in $H$.
\end{enumerate}
Moreover, if $H$ is not $11k\r'$-extremal, then $H$ contains a perfect matching if and only if $(\cP, L_{\cP}^{\mu}(H))$ is soluble.
\end{theorem}

We use following lemma which can be easily derived from Theorem~\ref{lem:PL}.

\begin{lemma}\label{prop:estimate}
Suppose $k\ge 3$, $2\le p\le k$ and $1/n \ll \beta, \mu_0 \ll \r' \ll 1/k$.
Let $H$ be a $k$-graph such that $\delta_{k-1}(H)\ge (1/p+11k\r')n$.
Then there exist $\mu_0^2 \le \mu\le \mu_0$ and a full pair $(\cP, L_{\cP}^{\mu}(H))$ with $\cP=\{V_1, \dots, V_d \}$ such that
\begin{enumerate}[\emph{(}i\emph{)}]
\item $d< p$, $|V_i| > n/p$ for any $i\in [d]$, and 
\item $H$ contains a perfect matching if and only if $(\cP, L_{\cP}^{\mu}(H))$ is soluble.
\end{enumerate}
\end{lemma}

\begin{proof}
Let $\mu$ and the full pair $(\cP, L_{\cP}^{\mu}(H))$ be returned by Theorem~\ref{lem:PL}.
First assume that $H$ is $11k\r'$-extremal, namely, $H$ contains an independent set of size at least $(1-11k\r')\frac{k-1}k n$.
Thus any $(k-1)$-tuple $S$ of vertices in this independent set has degree at most $n - (1-11k\r')\frac{k-1}k n = \frac nk + 11(k-1)\r' n < \delta_{k-1}(H)$, a contradiction.
Thus ($ii$) holds.

Now assume that there exists $i\in [d]$ such that $|V_i| \le n/p$. Write $L:=L_{\cP}^{\mu}(H)$.
Note that the number of edges that contain at least $k-1$ vertices in $V_i$ is at least $\frac1k\binom{|V_i|}{k-1}\delta_{k-1}(H) > d \mu n^k$, where we used $|V_i|\ge n/(2k)$ by (ii) and $\mu \ll 1/k$.
Thus there exists $j\in [d]$ such that $(k-1)\bfu_i+\bfu_j\in L$.
Let $x$ denote the number of $(k-1)$-sets $S$ such that $\bfi_{\cP}(S)=(k-2)\bfu_i+\bfu_j$.
By the minimum codegree condition and $\mu \ll \r'$, we have
\[
\sum_{\bfi_{\cP}(S)=(k-2)\bfu_i+\bfu_j} |N(S)\setminus V_i| \ge \sum_{\bfi_{\cP}(S)=(k-2)\bfu_i+\bfu_j} (\deg(S) - |V_i|) \ge x (1/p+11k\r')n - x |V_i| > d^k k \mu n^k.
\]
Since there are at most $d^k$ $k$-vectors and each edge is counted at most $k$ times in the above inequality, we conclude that there exists $i'\neq i$ such that there are at least $\mu n^k$ edges of $H$ with index vector $(k-2)\bfu_i+\bfu_j + \bfu_{i'}$, namely, $(k-2)\bfu_i+\bfu_j + \bfu_{i'}\in L$.
Thus $\bfu_i - \bfu_{i'}\in L$, contradicting that $L$ is transferral-free.
So we conclude that for any $i\in [d]$, $|V_i| > n/p$, which implies $d<p$ immediately.
So $(i)$ holds.
\end{proof}

Now we are ready to prove Theorems~\ref{thm:main1} and~\ref{thm:main2}.
Note that it remains to show that $(\cP, L_{\cP}^{\mu}(H))$ is soluble.

\begin{proof}[Proof of Theorem~\ref{thm:main1}]
Given an integer $k\ge 3$, let $p$ be the smallest prime factor of $k$.
By the monotonicity, we may assume that $\r\ll 1/k$.
Let $\r'=\r/11k$.
Moreover, suppose that $0<1/n\ll \beta, \mu_0 \ll \r \ll 1/k$ and $n\in k\mathbb{N}$.
Let $H=(V, E)$ be an $n$-vertex $k$-graph such that $\alpha(H)\le n/p$ and $\delta_{k-1}(H)\ge (1/p+\r) n$. 
Let $\mu$ and the full pair $(\cP, L_{\cP}^{\mu}(H))$ be returned by Lemma~\ref{prop:estimate}, where $\cP=\{V_1, \dots, V_d \}$.
So $d< p$ and $|V_i| > n/p$ for any $i\in [d]$.
Note that since $n\in k\mathbb{N}$, $\bfi_\cP(V)\in L_{\max}^d$.

It remains to show that $(\cP, L)$ is soluble, where $L:=L_{\cP}^\mu(H)$.
Consider the following set of vectors $\{\bfi_\cP(V)$, $k\bfu_1$, \dots, $k\bfu_d\}$.
Since by Lemma~\ref{lem:group}, $|G(L)|=d$, thus the set of vectors above contains two vectors that lie in the same coset of $L$, i.e., their difference is in $L$.
First assume that for some $i\in [d]$, $\bfi_\cP(V) - k\bfu_i\in L$.
Because $\alpha(H)\le n/p< |V_i|$, $H[V_i]$ contains an edge $e$.
Let $M=\{e\}$ be the matching of this single edge, then we have $\bfi_{\cP}(V\setminus V(M)) = \bfi_\cP(V) - k\bfu_i\in L$ and thus $(\cP, L)$ is soluble.
So we may assume that there exists $i, j\in [d]$, $i\neq j$, such that $k(\bfu_i - \bfu_j) \in L$.
In this case, consider the following set of vectors $\{\mathbf{0}$, $\bfu_i - \bfu_j$, $2(\bfu_i - \bfu_j), \dots, d(\bfu_i - \bfu_j)\}$, where $\mathbf{0}$ denotes the all $0$ vector.
Similarly by $|G(L)|=d$, some two of the vectors above lie in the same coset of $L$ and thus there exists $r\in [d]$ such that $r(\bfu_i - \bfu_j) \in L$.
Since $r\le d<p$, $r$ and $k$ are co-prime, namely, there exist integers $a, b$ such that $ar+bk=1$.
This implies that $a r(\bfu_i - \bfu_j) + b k(\bfu_i - \bfu_j) = \bfu_i - \bfu_j \in L$, contradicting that $L$ is transferral-free.
\end{proof}

Throughout the rest of the paper, for distinct vertices $x, y$ and a vertex set $U$, we write $\deg(xy, U) = |N(\{x, y\})\cap U|$.


\begin{proof}[Proof of Theorem~\ref{thm:main2}]
By the monotonicity, we may assume that $\r\ll 1$.
Let $\r'=\r/33$.
Moreover, suppose $0<1/n\ll \beta, \mu_0 \ll \r \ll 1$ and $n\in 3\mathbb{N}$.
Let $H=(V, E)$ be an $n$-vertex $3$-graph with $\delta_{2}(H)\ge (1/3+\r) n$ which is induced $K_4^-$-free. 
Let $\mu$ and the full pair $(\cP, L_{\cP}^{\mu}(H))$ be returned by Lemma~\ref{prop:estimate} with $k=p=3$, where $\cP=\{V_1, \dots, V_d \}$.
So $d\in \{1, 2\}$ and $|V_i| > n/3$ for any $i\in [d]$.
Note that since $n\in 3\mathbb{N}$, $\bfi_\cP(V)\in L_{\max}^d$.

It remains to show that $(\cP, L)$ is soluble, where $L:=L_{\cP}^\mu(H)$.
Note that if $d=1$, then by taking $M=\emptyset$ we see $(\cP, L)$ is soluble.
So assume $d=2$ and $\bfi_\cP(V)\notin L$ (otherwise we are done with $M=\emptyset$).
Because $(\cP, L)$ is full, it is easy to see that $(0,3), (2,1)\in L$ and $(3, 0), (1,2)\notin L$ (up to the permutation of the parts of $\cP$).
Since $|G(L)|=2$, the vectors $(3, 0), (1,2)$ and $\bfi_\cP(V)$ belong to the same coset of $L$.
Thus if $H$ contains any edge $e$ with index vector $(3,0)$ or $(1,2)$, then $\bfi_{\cP}(V\setminus e)\in L$ and we are done.
So we may assume that no such edge exists, and in particular, $V_1$ is independent.

Note that $|V_1|< 2n/3$ because $|V_2|> n/3$.
Thus we have $\delta_2(H)\ge (1/3+\r)n > |V_1|/2$.
Now fix any vertex $v\in V_2$.
For any vertex $v'\in V_1$, we know that $\deg(v v', V_1)= \deg(v v')\ge \delta_2(H)> |V_1|/2$.
Thus we have $|N_H(v)\cap \binom{V_1}2|>|V_1|^2/4$.
So by Mantel's Theorem, as a graph on $V_1$, $N_H(v)\cap \binom{V_1}2$ contains a triangle $v_1, v_2, v_3$.
Moreover, since $V_1$ is independent, $H[v, v_1, v_2, v_3]$ forms an induced copy of $K_4^-$, which is a contradiction.
So the proof is completed.
\end{proof}

\section{Almost perfect matchings and proofs of Theorems~\ref{thm:NPM2} and~\ref{thm:main3}}

In this section we first prove the following theorems on almost perfect matchings and tilings.

\begin{theorem}\label{thm:apm}
For any $\r>0$, the following holds for all sufficiently large $n$.
Let $H$ be an $n$-vertex $3$-graph such that $H$ is induced $K_4^-$-free and $\delta_{2}(H) \ge (1/6+\r)n$.
Then $H$ contains a matching that covers all but at most $5/\r$ vertices. 
\end{theorem}

\begin{proof}
Clearly we can assume that $\r\le 5/6$.
Assume that $M=\{e_1,\dots, e_m\}$ is a maximum matching in $H$ and let $U=V(H)\setminus V(M)$.
First observe that for any $v\in V(M)$, $\deg(v, U)\le |U|^2/4$.
Indeed, otherwise, by Mantel's Theorem $N_H(v)\cap \binom{U}2$ (as a graph) contains a triangle on $u_1, u_2, u_3$.
On the other hand by the maximality of $M$, $U$ is an independent set.
Thus $H[\{v, u_1, u_2, u_3\}]$ is an induced copy of $K_4^-$, a contradiction.

Let $D$ be the set of vertices $v\in V(M)$ such that $\deg(v, U)> 2|U|$ and observe that $|D|\le m$.
Indeed, we show $|D\cap e_i|\le 1$ for $i\in [m]$.
If $v', v''\in D\cap e_i$ for some $i\in [m]$, then by the definition of $D$, we can greedily find two disjoint edges $\{v', u_1, u_2\}, \{v'', u_3, u_4\}$ in $H$ such that $u_1,\dots, u_4 \in U$.
By replacing $e_i$ in $M$ by these two edges, we get a matching in $H$ larger than $M$, a contradiction.

Let $e_H(V(M), U, U)$ denote the number of edges of $H$ that contain one vertex in $V(M)$ and two vertices in $U$, then we have
\[
\binom{|U|}{2} (1/6+\r)n \le \binom{|U|}2 \delta_2(H) \le e_H(V(M), U, U)\le |D|\cdot \frac{|U|^2}{4} + 3m \cdot 2 |U|.
\]
Applying $|D|\le m$ and $n\ge 3m$ and then dividing both sides by $m|U|/4$, we get $(|U|-1) (1+6\r) \le |U| +24$, which implies that $|U|\le 25/(6\r)+1\le 5/\r$.
\end{proof}

We use the following result due to Goodman~\cite{Goodman} and Moon and Moser~\cite{MoMo62} on the number of triangles in a graph.

\begin{theorem}\cite{Goodman, MoMo62}\label{thm:Goodman}
A graph with $n$ vertices and $m$ edges contains at least $m(4m-n^2)/3n$ triangles.
\end{theorem}

\begin{theorem}\label{thm:apm2}
For any $\r>0$, the following holds for all sufficiently large $n$.
Let $H$ be an $n$-vertex $3$-graph such that $H$ is induced $K_4^-$-free and $\delta_{2}(H) \ge (1/8+\r)n$.
Then $H$ contains a $Y$-tiling that covers all but at most $40/\r+1$ vertices. 
\end{theorem}

\begin{proof}
Assume that $M=\{F_1,\dots, F_m\}$ is a maximum $Y$-tiling in $H$ and let $U=V(H)\setminus V(M)$.
Assume that $|U| > 40/\r+1$.
First observe that for any $v\in V(M)$, $\deg(v, U)\le (1/2+\r)\binom{|U|}2$.
Indeed, otherwise, by Theorem~\ref{thm:Goodman}, $N_H(v)\cap \binom{U}2$ (as a graph on $U$) contains 
\[
(1/2+\r)\binom{|U|}{2} \left(4(1/2+\r)\binom{|U|}{2} - |U|^2 \right)\frac{1}{3|U|} \ge \frac{1}{6}\binom{|U|}{2}\cdot (2\r |U| - (1+2\r)) \ge \frac{\r}2 \binom{|U|}3
\]
triangles, where we used $\r |U| \ge 1$ in the last inequality.
Since $H$ is induced $K_4^-$-free, each such triangle $u_1 u_2 u_3$ must form an edge in $H$.
By $|U|> 6/\r+2$, we see that $|E(H[U])|\ge \frac{\r}2 \binom{|U|}3 > \binom{|U|}2$, which, by averaging, implies that there is a pair of vertices with degree at least $2$ in $U$, i.e., a copy of $Y$ in $U$, contradicting the maximality of $M$.

Let $D$ be the set of vertices $v\in V(M)$ such that $\deg(v, U)> 5|U|$ and observe that $|D|\le m$.
Indeed, we show $|D\cap F_i|\le 1$ for $i\in [m]$.
If $v', v''\in D\cap F_i$ for some $i\in [m]$, then by $\deg(v', U)> 5|U|$, there is a vertex $u_1$ such that $N(v' u_1, U)\ge 2$, i.e., there exists $u_2, u_3$ such that $H[v', u_1, u_2, u_3]$ forms a copy of $Y$.
Moreover, note that $\deg(v'', U\setminus \{u_1, u_2, u_3\}) > 2 |U|$, and similarly we can pick vertices $u_4, u_5, u_6\in U\setminus \{u_1, u_2, u_3\}$ such that $H[v'', u_4, u_5, u_6]$ forms a copy of $Y$.
By replacing $F_i$ in $M$ by these two copies of $Y$, we get a $Y$-tiling in $H$ larger than $M$, a contradiction.

Moreover, note that by the maximality of $M$, for any pair of vertices $u v$ in $U$, we have $\deg(u v, U)\le 1$.
Let $e_H(V(M), U, U)$ be the number of edges of $H$ that contain one vertex in $V(M)$ and two vertices in $U$, then we get
\[
\binom{|U|}{2} (1/8+\r/2)n \le \binom{|U|}2 (\delta_2(H)-1)\le e_H(V(M), U, U)\le |D|\cdot (1/2+\r)\binom{|U|}2 + 4m \cdot 5 |U|.
\]
Applying $|D|\le m$ and $n\ge 4m$ and then dividing both sides by $m|U|/4$, we get $(|U|-1) (1+4\r) \le (|U|-1)(1+2\r) +80$, which implies that $|U|\le 40/\r+1$, a contradiction.
\end{proof}

Theorem~\ref{thm:NPM2} can be proved by combining Theorem~\ref{thm:apm} and Lemma~\ref{lem:abs2} under the same fashion as the proof of Theorem~\ref{thm:NPM1}, and thus we omit the details.
To prove Theorem~\ref{thm:main3}, we need the following absorbing lemma proved by Czygrinow, DeBiasio and Nagle~\cite{CDN}.

\begin{lemma}\cite[Lemma 3.3]{CDN}
For all $\r, \delta>0$, there exists $\alpha>0$ so that the following holds for all sufficiently large $n\in 4\mathbb{N}$.
Let $H$ be an $n$-vertex $3$-graph with $\delta_2(H)\ge \delta n$.
Then there exists $A\subseteq V(H)$ of size $|A|\le \r n$ so that for every $W\subseteq V(H)\setminus A$ of size $|W|\le \alpha n$ for which $|A\cup W|$ is divisible by $4$, $H[A\cup W]$ contains a perfect $Y$-tiling.
\end{lemma}

By Theorem~\ref{thm:apm2} and the lemma above, the proof of Theorem~\ref{thm:main3} is a standard application of the absorbing method and is thus omitted.

\section{Sharpness of the results}

In this section we construct $k$-graphs that show the sharpness of our results.

\subsection{Sharpness of Theorems~\ref{thm:main1} and~\ref{thm:main2}}
We use a construction of Mycroft~\cite{My14}.
For any integer $k\ge 3$, let $p$ be the smallest prime factor of $k$.
We first define the following set of vectors in $\mathbb{Z}_p^p$.
For $1\le i<p$, let $\bfv_i=\bfu_i + (i-1)\bfu_p$.
Let $L$ be the set of vectors generated by $\bfv_1,\dots, \bfv_{p-1}$.
The following claim was proved in~\cite{My14} and we include a proof for completeness.

\begin{claim}\label{claim:L}
$L$ is transferral-free. Moreover, for any $(k-1)$-vector $\bfv$, there exists $i\in [p]$ such that $\bfv + \bfu_i\in L$.
\end{claim}

\begin{proof}
We first show that $L$ is transferral-free.
Write $\bfv_p = \bfu_p + (p-1)\bfu_p = \mathbf{0}\in L$.
Now assume that there exist distinct $i,j\in [p]$ such that $\bfu_i - \bfu_j\in L$.
Since $\bfv_i, \bfv_j\in L$, we have
$\bfv_i - \bfv_j = \bfu_i - \bfu_j + (i-j) \bfu_p \in L$.
Thus we conclude that $(i-j) \bfu_p \in L$.
This is a contradiction because when $i- j\not\equiv 0 \pmod p$, the equation
\[
(i-j) \bfu_p \equiv \sum_{1\le i'<p} a_{i'} \bfv_{i'} \, \pmod p
\]
has no integer solution.

Next fix a $(k-1)$-vector $\bfv=(a_1,\dots, a_p)$, and let $b:=\sum_{j\in [p-1]} a_j(j-1)$.
Note that
\[
\bfv + (b-a_p) \bfu_p = (a_1, \dots, a_{p-1}, b) = \sum_{j\in [p-1]} a_j \bfv_j \in L.
\]
Let $i\in [p]$ such that $b-a_p\equiv -(i-1) \pmod p$ and note that $\bfv + \bfu_i = (\bfv - (i-1)\bfu_p) + \bfv_i\in L$.
\end{proof}

Let $\cP=\{V_1, V_2, \dots, V_p\}$ be a partition of $V$ such that $|V_1|+|V_2|+\cdots +|V_p|=n$, $|V_i| = n/p \pm 1$ for $i\in [p]$ and $\bfi_\cP(V)\notin L$ modulo $p$ (on a co-ordinate by co-ordinate basis).
Note that the restriction on the sizes is possible because $L$ is transferral-free: either $(n/p, \dots, n/p)$ or $(n/p-1, n/p+1, n/p,\dots, n/p)$ is not in $L$ modulo $p$.
Let $H_k$ be the $k$-graph with vertex set $V$ whose edges are $k$-tuples $e$ such that $\bfi_{\cP}(e)\in L$ modulo $p$.
Observe that $\delta_{k-1}(H_k) \ge n/p-k$.
Indeed, for any $(k-1)$-set of vertices $S$, by the second part of Claim~\ref{claim:L}, there exists $i\in [p]$ such that $\bfi_{\cP}(S) + \bfu_i\in L$.
By the definition of $H_k$, we have $V_i\setminus S\subseteq N_{H_k}(S)$ and thus $\deg_{H_k}(S)\ge |V_i\setminus S| \ge n/p-1-(k-1)=n/p-k$.
Moreover, for any matching $M$ in $H_k$, we have $\bfi_\cP(V(M))=\sum_{e\in M}\bfi_\cP(e)\in L$ modulo $p$.
Because $\bfi_\cP(V)\notin L$, $V(M)\neq V$ and thus $M$ is not perfect.
Moreover, since $p\mid k$, for $i\in [p]$, all the $k$-tuples $S$ in $V_i$ are edges of $H_k$ (because $\bfi_{\cP}(S)\equiv \mathbf{0}\in L$ modulo $p$), so $\a(H_k) < pk$.
At last, when $k=p=3$, $L$ is generated by $(1,0,0)$ and $(0,1,1)$ (modulo $3$).
So
\[
E(H_3) = \left\{e\in \binom V3:\bfi_{\cP}(e)\in\{(3,0,0), (0,3,0), (0,0,3), (1,1,1)\}\right\}, 
\]
and it is straightforward to check that $H_3$ is induced $K_4^-$-free.
This shows the sharpness of the minimum codegree conditions in Theorems~\ref{thm:main1} and~\ref{thm:main2}.

\subsection{Sharpness of Theorems~\ref{thm:NPM2} and~\ref{thm:main3}}
For the sharpness of Theorems~\ref{thm:NPM2} and~\ref{thm:main3}, consider the following $3$-graph.
Let $m,n$ be integers such that $m\le n/2$.
Let $A$ and $B$ be two disjoint vertex sets such that $|A|=m$ and $|B|= n - |A|$.
Let $E$ be the set of the triples with one vertex in $A$ and two vertices in $B$ generated as follows:
for each vertex $a$ in $A$, pick a random bipartition $B_1\cup B_2$ of $B$ such that $||B_1| - |B_2||\le 1$, and add all triples $\{a b_1 b_2: b_1\in B_1, b_2\in B_2\}$ to $E$.
Now let $H(m, n)$ be the $3$-graph with vertex set $A\cup B$ and whose edges are all triples in $E$ and all triples with three vertices in $A$.
Note that by standard concentration results (e.g., Chernoff's bound), if $m=\Theta(n)$, then for any pair of distinct vertices $b, b'\in B$, it holds that $\deg_{H(m, n)}(b b') = m/2 - o(n)$.
Then it is straightforward to check that $\delta_2(H(m, n)) = m/2 - o(n)$.
Finally, we show that $H(m, n)$ is induced $K_4^-$-free.
Indeed, let $F$ be an induced copy of $K_4^-$ in $H(m, n)$ with vertex set $W$.
Since no edge of $H(m, n)$ contains exactly two vertices in $A$, then $|W\cap A|\notin \{2,3\}$.
Moreover, since $H(m, n)[A]$ is complete and $H(m, n)[B]$ is independent, the only possible case is $|W\cap A| = 1$.
However, by the definition of $E$, for any $v\in A$, $N(v)\cap \binom{B}{2}$ is triangle-free, and since $H(m, n)[B]$ is independent, there are at most two edges on $W$.
So we get a contradiction and we are done.

Taking $m=n/3 - 2$, since $|A| = n/3-2$ and $B$ is independent, the largest matching in $H(n/3-2, n)$ covers $n-6$ vertices.
Since $\delta_2(H(n/3-2, n)) = (1/6 - o(1))n$, $H(n/3-2, n)$ shows the sharpness of the minimum codegree condition in Theorem~\ref{thm:NPM2}.
Similarly, taking $m=n/4-1$ shows the sharpness of the minimum codegree condition in Theorem~\ref{thm:main3}.

\section{Concluding Remarks}
In this paper we studied Dirac-type problems under additional structural assumptions.
Let us mention some remarks and further directions.

\subsection{$k$-graphs with small independence number or uniformly dense}
The study of (hyper)graphs with sublinear independence number actually dates back to the classic Ramsey-Tur\'an theory.
The original motivation is that the extremal examples for Tur\'an problems are well-structured, which is `as far as possible' from the `randomlike' structures which arise in bounds on Ramsey numbers.
From this point of view, besides the Ramsey-Tur\'an problems, Erd\H{o}s and S\'os~\cite{ErSo82} also posed some questions about $k$-graphs with uniform edge density, namely, for the $n$-vertex $k$-graphs $H$ such that

\begin{itemize}
\item[$(\dagger)$] for every subset $W$ of $V(H)$ of size at least $\e n$, $H[W]$ contains at least $c |W|^k$ edges,
\end{itemize}
where $1/n\ll \e \ll c$.
We remark that the minimum codegree condition in Theorem~\ref{thm:main1} cannot be improved even by replacing the independence number condition with the stronger condition~$(\dagger)$, because the lower bound construction presented in Section 4 satisfies~$(\dagger)$ with $c=1/k^k$.

It is natural to investigate the perfect $F$-tiling problem in $k$-graphs with small independence number or those satisfying~$(\dagger)$.
It is also natural to start with the case $k=3$ and some small $F$.
Let $K_4^3$ be the unique $3$-graph with $4$ vertices and $4$ edges (so the complete $3$-graph on $4$ vertices).

\begin{problem}\label{p1}
For $F=K_4^3$ or $K_4^-$, determine the minimum codegree threshold for perfect $F$-tilings in $3$-graphs $H$ with $\alpha(H)=o(n)$ or $3$-graphs $H$ satisfying $(\dagger)$ (for some nontrivial $c$).
\end{problem}

Note that for general $3$-graphs, Keevash and Mycroft~\cite{KM1}, and independently, Lo and Markstr\"om~\cite{LM1} showed that the threshold for perfect $K_4^3$-tiling is $(3/4+o(1))n$ (in~\cite{KM1} the exact threshold was given).
For perfect $K_4^-$-tilings, Lo and Markstr\"om~\cite{LM2} showed that the threshold is $(1/2+o(1))n$ (the exact threshold was determined recently in~\cite{HLTZ_K4}).
However, the examples supporting the lower bounds in both results contain large independent sets.
Here we present some nontrivial lower bound constructions satisfying $(\dagger)$ (thus also having small independence numbers), essentially from~\cite{Pik}.

It is known that there is a $K_4^3$-free $3$-graph $H_1(n)$ on $n$ vertices with minimum codegree $(1/2 - o(1))n$, and there is a $K_4^-$-free $3$-graph $H_2(n)$ on $n$ vertices with minimum codegree $(1/4 - o(1))n$, and both $H_1(n)$ and $H_2(n)$ are `randomlike' and hence satisfy~$(\dagger)$ with $c=1/2-o(1)$ and $c=1/4-o(1)$, respectively.
Indeed, let $T$ be a random tournament on $n$ vertices. For $H_1(n)$, fix a linear order of the vertices and add a triple $\{x, y, z\}$ with $x<y<z$ to $H_1(n)$ if out of the two pairs $\{x, y\}$ and $\{x, z\}$ exactly one is directed toward $x$ in $T$; for $H_2(n)$, let $H_2(n)$ be the set of triples that span a cyclic triangle in $T$. The estimates on minimum codegrees follow from standard concentration results and we omit them.

Now let $V=A\cup B$ of size $n\in 4\mathbb N$ such that $A\cap B=\emptyset$, $|A|=n/4-1$.
Let $E_A$ be the set of the triples in $V$ that intersect $A$.
Let $H_1'$ be the $3$-graph on $V$ with edge set $E_A\cup H_1(|B|)$, where $H_1(|B|)$ is defined on $B$.
Note that $\delta_2(H_1')\ge (5/8+o(1))n$ and $H_1'$ satisfies~$(\dagger)$ with $c=1/2-o(1)$.
To see that $H_1'$ has no perfect $K_4^3$-tiling, note that since $H_1'[B]$ is $K_4^3$-free, every copy of $K_4^3$ in $H_1'$ must contain at least one vertex in $A$, but $|A|<n/4$.
Let $H_2'$ be the $3$-graph on $V$ with edge set $E_A\cup H_2(|B|)$, where $H_2(|B|)$ is defined on $B$.
Similarly $\delta_2(H_2')\ge (7/16+o(1))n$, $H_2'$ satisfies~$(\dagger)$ with $c=1/4-o(1)$, and it has no perfect $K_4^-$-tiling.

\subsection{$k$-graphs without induced copy of a certain $k$-graph}

In this paper we studied induced $K_4^-$-free $3$-graphs. It is natural to investigate the problems for $k$-graphs with $k\ge 3$ and other forbidden subgraphs.

Let us briefly mention the case for $3$-graphs which are $K_4^3$-free (we omit the word `induced' here because every copy of $K_4^3$ must be induced).
Consider the following example.
Let $V=V_1\cup V_2$ of size $n\in 3\mathbb N$ such that $V_1\cap V_2=\emptyset$, $|V_2|\in [n/3, n/3+1]$ is odd.
Let $E_1$ be the set of the triples in $V$ that contain exactly one vertex from $V_1$ and two vertices from $V_2$.
Let $H'$ be the $3$-graph on $V$ with edge set $E_1\cup H_1(|V_1|)$ on $V_1$, where $H_1(|V_1|)$ is the $3$-graph mentioned above.
It is easy to see that $H'$ contains no copy of $K_4^3$, and it contains no perfect matching (because each edge of $H'$ contains $0$ or $2$ vertices in $V_2$ but $|V_2|$ is odd).
Moreover, $\delta_2(H')=(1/3-o(1))n$.

Here is another example which achieves a better minimum codegree.
To be brief, we only show it for $n\in 3\mathbb N\setminus 9\mathbb N$.
Let $V_0, V_1, V_2$ be three disjoint vertex sets each of size $m\notin 3\mathbb N$.
let $V=V_0\cup V_1\cup V_2$ and $n=3m$.
For $i=0,1,2$, let $E_i$ be the set of the triples that contain two vertices from $V_i$ and one vertex from $V_{i+1}$, where we treat $V_3=V_0$.
Then let $H''$ be the $n$-vertex $3$-graph on $V$ with the edge set $E_0\cup E_1\cup E_2$.
Clearly $\delta_2(H'')=m-1=n/3-1$ and $H''$ contains no copy of $K_4^3$.
Assume that $H''$ contains a perfect matching, which uses $a_i$ edges from $E_i$, $i=0,1,2$.
Thus we have $2a_0+a_2=2a_1+a_0=2a_2+a_1=m$, which implies that $a_0=a_1=a_2=m/3\notin \mathbb N$, a contradiction.

\begin{conjecture}\label{conjk43}
For any $\r>0$, the following holds for sufficiently large integer $n\in 3\mathbb{N}$.
Suppose $H$ is an $n$-vertex $3$-graph such that $H$ does not contain any copy of $K_4^3$ and $\delta_{2}(H) \ge (1/3+\r)n$, then $H$ contains a perfect matching.
\end{conjecture}

\section*{Acknowledgement}
The author is indebted to Hi\d{\^{e}}p H\`an, Yoshiharu Kohayakawa and Liming Xiong for helpful discussions.
Part of this research was carried out while the author was visiting Universidad de Santiago de Chile, and the author would like to thank the university for the nice working environment.
The author also would like to thank two anonymous referees for the valuable comments that improve the presentation of the paper.


\end{document}